\documentclass{amsart}
\usepackage[mathscr]{eucal}
\usepackage{amssymb, amsmath,array, amscd}
\usepackage{enumerate}

\newtheorem{thm}{Theorem}[section]

\newtheorem{prop}[thm]{Proposition}
\newtheorem{lemma}[thm]{Lemma}

\newtheorem{remark}[thm]{Remark}
\newtheorem{cor}[thm]{Corollary}

\newtheorem{example}[thm]{Example}

\begin{document}

\title{A new bound on the number of special fibers in a pencil of curves}

\author{S. Yuzvinsky}
\address{Departament of Mathematics \\ University of Oregon \\ Eugene \\ OR 94703 USA}
\email { yuz@uoregon.edu }

\begin{abstract}
In paper \cite{PY} it was proved 
that any pencil of plane curves of degree $d>1$
with irreducible generic fiber can have at most five completely reducible fibers although no examples with five such fibers were ever
found. Recently Janis Stipins has proved in \cite{St} that  if a pencil has the base of
$d^2$ points then it cannot have five completely reducible fibers. In this paper
we generalize Stipins' result to arbitrary pencils. We also include into consideration more general special fibers that are the unions of lines and non-reduced
curves. These fibers are important for characteristic varieties of line complements.
 \end{abstract}

\maketitle

\date{\today}

\section{Introduction}
The main object of this note is a 1-dimensional system 
of homogeneous polynomials
of three variables of some degree $d>1$ or the respective projective system (pencil) 
of curves in $\mathbb P^2$. We identify these systems and use the symbol $\mathcal P$
for both. If the generic member (fiber) of $\mathcal P$ is irreducible then the set of
reducible fibers is finite. There are various known
upper bounds on the size of this set
that depend on $d$, see, for instance, \cite{Vi,Bo}. 

We concerned mainly with \emph{completely reducible fibers}, i.e., unions of lines not necessarily
reduced. Certain results for these fibers with some extra conditions have appeared in 
classical mathematical literature (see \cite{Ha,Hal}). New interest to completely reducible 
fibers was ignited by papers \cite{LY,FY} where the union of these fibers 
 was characterized by topological 
property of its complement and by combinatorics of intersection of lines.
The best upper bound on the number $k$ of the completely reducible fibers 
 was found in \cite{PY} where the inequality $k\leq 5$ was proved without any extra conditions
 on $\mathcal P$. There were however no examples known with $k=5$.
 
 Then a better bound appeared in dissertation \cite{St} where the inequality $k\leq 4$
 was proved assuming that any two different completely reducible fibers
 intersect transversally. The condition is equivalent  to all the lines of those
 fibers forming a $(k,d)$-net (see \cite{Yu}). 
 
 The goal of this note is to generalize
 the method from \cite{St} in order to prove that $k\leq 4$ without extra conditions.
 Also a more general
  kind of fibers appeared in studies by A.Dimca of characteristic varieties 
 for line arrangement complements (see \cite{Di}). These are fibers that are unions
 of lines and non-reduced curves. We called these fibers \emph {special}
 and  include them into considerations.
 
 Let us remark that there are many examples of pencils with $k=3$
 (for instance, see \cite{FY,PY,St}) but the only one (up to isomorphism)
 with $k=4$. This is the well-known Hesse pencil generated by 
 a smooth plane cubic and its Hessian. We consider the case $k=4$
 in Theorem \ref{mainer} and Problem at the end of the note. 
 
  Since the dissertation \cite{St} is not published we do not assume that its proofs are
 available to the reader.  As the result this note has proofs of three different types.
 Some of them, like Lemma \ref{dolga} (attributed by Stipins to Dolgachev) and
 Proposition \ref{stipins}, are essentially taken from \cite{St}, some, like Theorem \ref{noether},
 made shorter and fit the more general situation, others, like Theorem \ref{mainer},
 are new.
 
 The author is grateful to J.Stipins for providing his dissertation and to
I.Dolgachev for useful discussions.

\section{Hessians of a pencil of curves}
Recall from Introduction that as a rule we do not distinguish between a plane curve
and its defining polynomial. Let $\mathcal P$ be a pencil of curves of degree $d$ ($d>1$) whose
generic fiber is irreducible. We assume immediately that $\mathcal P$ has at least four special
fibers $F_i$, $i=1,2,3,4$, three of which, say $F_i,\ i=1,2,3$ are completely reduclible.
Clearly any two fibers cannot have a common component whence the base $B$ of $\mathcal P$ (i.e., the mutual intersection of all the fibers) is finite. Also it follows from the conditions on $\mathcal P$ that
$|B|>1$, i.e., $\mathcal P$ does not lie in a pencil of lines.
Choosing the fibers $F_1$ and $F_2$ as generators we have
$F_i=\lambda_iF_1+\mu_iF_2$, $i=3,4$ for some $[\lambda_i:\mu_i]\in \mathbb P^1$.
 For symmetry we put $(\lambda_1,\mu_1)=(1,0)$ and $(\lambda_2,\mu_2)=(0,1)$. 
We also put $F(\lambda,\mu)=\lambda F_1+\mu F_2$ for all $\lambda$ and $\mu$.

Now we study the family parametrized by $[\lambda:\mu]\in\mathbb P^1$ of the Hessians of curves
$F(\lambda,\mu)$. Let us recall that the Hessian of a plane curve (or a homogeneous polynomial of three variables $x, y,z$) is the polynomial $H(F)$ that is the determinant of the
$3\times 3$-matrix of the second derivatives of $F$. The points of $F\cap H(F)$ are precisely the points of $F$ that are either multiple or flexes (see, for instance \cite{Fu}). In particular $H(F)$ contains $F$ if and only if
$F$ is a special fiber.
Also $H(F)$ is identically zero
if and only if $F$ is the union of concurrent lines. We put
$H(\lambda,\mu)=H(F(\lambda,\mu))$. Notice that $H(\lambda,\mu)$ is a homogeneous
polynomial of degree $3d-6$ in $x,y,z$ whose coefficients are cubic forms in $\lambda$ and $\mu$.

\begin{lemma}
\label{dolga}
For every $[\lambda:\mu]\in\mathbb P^1$ the curve $H(\lambda,\mu)$ passes through every
$P\in B$.
\end{lemma}
\begin{proof} For a fixed point $P\in B$ the polynomial $H(\lambda,\mu)$ evaluated at $P$ is a cubic form in $\lambda$ and $\mu$ that has at least 4 zeros at $(\lambda_i,\mu_i)$. Thus this form is identically 0 which implies the statement.
\end{proof}

Lemma \ref{dolga} is false for a pencil with just three completely reducible fibers.
\begin{example}
Consider the pencil $\{\lambda x^2(y^2-z^2)+\mu y^2(x^2-z^2)\}$ (closely related to the Coxeter group of type $B_3$). Then $H(1,1)$ is not 0 at the base point $[1:1:1]$. The pencil has three completely reducible fibers - the generators and $F(1,-1)=z^2(y^2-x^2)$.
\end{example}

In fact we need much stronger statement than Lemma \ref{dolga} - we need to check that 
the triple $(F_1,F_2,H(\lambda,\mu))$ satisfies at any $P\in B$ the conditions
 of the Max Noether Fundamental Theorem (\cite{Fu}, p. 120). Then we will be able
to apply that theorem. For that we need to discuss the multiplicities of the base 
points at the fibers of the pencil.

Fix $P\in B$ and without any loss assume that $P=[0:0:1]$.
Then the local ring $O_P(\mathbb P^2)$ can be identified with $\mathbb C_{(x,y)}[x,y]$. Also for any homogeneous polynomial $F$ of variables $x,y,z$ put $F^*=F(x,y,1)$ and write $F^*=\sum_{i\geq m}F^{(i)}$
where $F^{(i)}$ is the homogeneous component of degree $i$ of $F^*$ and
$F^{(m)}\not=0$. Here the non-negative
integer $m=m_P(F)={\rm mult}_P(F)$ is \emph{the multiplicity of $P$ on $F$}. 
The curve  $F^{(m)}$ in $\mathbb P^2$ is the union of 
concurrent lines with some multiplicities and each of them is a tangent line
 at $P$ to $F$.

In the rest of the paper any time we fix a point $P\in B$ we assume that $P=[0:0:1]$.

\begin{lemma}
\label{mult}

(i) The multiplicity of any $P\in B$ is the same on any fiber.

(ii) No two fibers have a common tangent line at any $P\in B$.

\end{lemma}

\begin{proof}
(i) Fix $P\in B$ and suppose there are fibers $J$ and $K$ such that $m=m_P(J)<m_P(K)$.
Then since $J$ and $K$ linearly generate all fibers we have $m_p(F)=m$ for every fiber $F$
different from $K$. Since  the number of the completely reducible fibers is at least 3
there exist two fibers, say $F_1$ and $F_2$,
among the completely reducible ones different from $K$ . Since for a completely reducible fiber $F$ the polynomial
$F^{(m)}$ is the product of all factors (with their multiplicities) of $F$ that pass through $P$ we conclude that $F^{(m)}_1$ and $F^{(m)}_2$ are linearly independent
over $\mathbb C$. Hence $F^{(m)}\not= 0$ for every fiber $F$ which contradicts the assumption about $K$.

(ii) Suppose there exist fibers $J$ and $K$ and a linear form $\alpha$ such that $\alpha$ divides
$J^{(m)}$ and $K^{(m)}$ where $m=m_P(F)$ for every fiber $F$. Then $\alpha$ divides
$F^{(m)}$ for every fiber $F$, in particular for completely reducible fibers $F_1$ and $F_2$. This implies
that $\alpha$ divides $F_1$ and $F_2$ which contradicts the assumption about $B$.
\end{proof}

\begin{remark}
Notice that the proof of Lemma \ref{mult} uses only that the pencil has three
completely reducible fibers. In a different form it was proved in
\cite{FY}.
\end{remark}

Now we study the multiplicity of the Hessians at the base points.
First we record a simple general lemma.
\begin{lemma}
\label{hessian}
 Let  $P=[0:0:1]$ be a point, $F$ a plane curve of degree $d$  that is not the union 
 of lines passing all through $P$,  and 
 $m=m_P(F)>1$. Let $H=H(F)$ be the Hessian of $F$. Then
 
  $m_P(H)=3m-4$ and $F^{(m)}$ divides $H^{(3m-4)}$.
\end{lemma}
\begin{proof}
It is clear from the definition that $m_P(H)\geq 3m-4$ and $H^{(3m-4)}z^{3(d-m)-2}=H(G)$
 where $G=F^{(m)}z^{d-m}$. 
 Notice that by condition on $F$ we have $d-m\geq 1$. Thus the curve $G$ is the union of non-concurrent lines whence it divides its Hessian which is not identically zero. Both statements follow.
\end{proof}
\begin{remark}
Notice that the previous lemma can be applied to every fiber $F$ of $\mathcal P$ and every $P\in B$ with $mult_P(F)>1$. Indeed if $F$
were the union of lines passing through $P$  then $P$ would have
been the only point in $B$ which is impossible by conditions on 
$\mathcal P$.
\end{remark}  
Now we can prove  what we need.
\begin{lemma}
\label{noether}
For all $[\lambda:\mu]\in\mathbb P^1$ Noether's conditions (see 
\cite{Fu}) are satisfied for $(F_1,F_2,H(\lambda,\mu))$ at every $P\in B$.
\end{lemma}
\begin{proof}
Fix $[\lambda:\mu]\in \mathbb P^1$ and consider several cases.

(i) $m_P(F_i)=1, i=1,2.$ The claim follows from Lemma \ref{dolga} and \cite{Fu}, Proposition 5.5.1(1).

(ii) $m=m_P(F_i)\geq 3, i=1,2.$ Then $m_P(F(\lambda,\mu))=m$ by Lemma \ref{mult}
whence by Lemma \ref{hessian} $m_P(H(\lambda,\mu))=3m-4\geq 2m-1$. The claim follows from
\cite{Fu}, Proposition 5.5.1(3).

(iii) $m_P(F_i)=2, i=1,2.$ 
By Lemma \ref{hessian} $m_P(H(\lambda,\mu))=2$ also.
Since no ready criterion applicable to this case is known to me we need to work with Noether's conditions directly.

To simplify notation put $F=F(\lambda,\mu)$ and $H=H(\lambda,\mu)$.
It suffices to prove that
$H^*\in (F_1^*,F_2^*)$ in $\mathbb C_{(x,y)}[x,y]$. 
 Using Lemma \ref{hessian} again we see that there exists $c\in\mathbb C^*$ such that
$H^{(2)}=c(\lambda F_1^{(2)}+\mu F_2^{(2)})$.
Now put $H_1=H-cFz^{2d-6}$ and notice that $m_P(H_1)\geq 3$.
Thus for $F_1,F_2,H_1$
the criterion in \cite{Fu}, Proposition 5.5.1(3) holds whence $H^*_1\in (F_1^*,F_2^*)$ 
in the local ring
$\mathbb C_{(x,y)}[x,y]$. Since this inclusion holds trivially for $H^*-H^*_1$ 
we complete the proof.
\end{proof}

\begin{remark} In this remark we briefly discuss the simple case $d=2$. First of all every point of the base
must have the multiplicity equal one on any fiber. Indeed if a point $P$ has multiplicity two  (the maximal possible) on $F_1$
then $P$ is the only point of $B$ which creates a contradiction (cf. the
previous remark). Thus $|B|=4$
which means that all the special fibers are completely reducible
and 
the union of all completely reducible fibers form a $(k,2)$-net for some $k>2$ 
(e.g., see \cite{FY}). 
It is an easy property of nets that the only $(k,2)$-net (up to linear
 isomorphism) is the $(3,2)$-net closely related to the Coxeter group
of type $A_3$ (for instance, see \cite{Yu}, p.1617). Because of that we will always assume in the rest of the paper that $d\geq 3$.
\end{remark}

Now we are ready to prove the main result about the family of the Hessians.
Recall that for each special fiber $F$ the Hessian $H(F)$ is divisible by $F$
and put $H(F_i)=C_iF_i$ for $i=1,2,3,4$ where each $C_i$ is a homogeneous polynomial of degree
$2d-6$.
\begin{thm}

\label{main}
There exist homogeneous polynomials $A$ and $B$ of degree $2d-6$ in $x,y,z$ whose coefficients are
homogeneous cubic polynomials in $\lambda$ and $\mu$ such that $H(\lambda,\mu)=AF_1+BF_2$ for all $[\lambda:\mu]\in\mathbb P^1$. 
Moreover these polynomials can be chosen so that

$${\rm(*) } A(\lambda_i,\mu_i)=\lambda_iC_i\  {\rm and}\  B(\lambda_i,\mu_i)=\mu_i C_i\ {\rm for}\  i=1,2,3,4.$$
\end{thm}

\begin{proof}
First of all Lemma \ref{noether} allows us to apply the Max Noether Fundamental Theorem which
implies the following. For every $(\lambda,\mu)\in\mathbb C^2$ there exist homogeneous polynomials $A=A(\lambda,\mu),B=B(\lambda,\mu)\in\mathbb C[x,y,z]$ of degree $2d-6$ such that $H(\lambda,\mu)=AF_1+BF_2$.

Now we prove the statement, namely that $A$ and $B$ can be chosen depending polynomially on $\lambda$ and $\mu$ and satisfying the conditions (*). For that we use elementary linear algebra over a polynomial ring. Indeed the coefficients of $A$ and $B$ are solutions of a finite system $\Sigma$ of linear equations whose coefficients
are complex numbers (in fact coefficients of polynomials
$F_i$, $i=1,2$) and the right-hand sides are values of cubic forms in $\lambda$ and $\mu$ (coefficients of 
polynomials $H(\lambda,\mu)$). It follows from the previous paragraph
 that this system
is consistent for each values of $\lambda$ and $\mu$. 
This implies that reducing the coefficient matrix
of $\Sigma$ to a row echelon form one can skip the zero rows and moving the free-variables columns
to the right-hand side obtain an equivalent system $\Sigma'$ with a
square non-degenerate matrtix
 of scalar coefficients.

Now we take care about the free variables. To get a solution 
polynomially depending on $\lambda$ and $\mu$
it is enough to choose the free variables depending on parameters as arbitrary homogeneous 
cubic polynomials. We want however to find a special solution satisfying the `boundary' conditions (*) at the
special values of parameters. Since the conditions (*) give a solution of the system $\Sigma'$
for these special values of parameters there are values of the free variables corresponding to them. Thus for each free
variable there exists a (unique) cubic form admitting each of these special values at the respective special value of parameters. Now combining cubic forms in the right-hand side and using Cramer's rule we obtain a solution with the required properties.
\end{proof}

\begin{remark} If $d\leq 5$ then the solution of the system $\Sigma$ is unique whence it is automatically polynomial and
satisfies the conditions (*). In this case there are no free variables. If however $d>5$ then one can easily find non-polynomial solutions by choosing free variables non-polynomially depneding on parameters.
\end{remark}

\begin{cor}
\label{form}
Theorem \ref{main} implies that
$$H(\lambda,\mu)=D\lambda F_1+E\mu F_2$$
where $D$ and $E$ are forms of degree $2d-6$ in $x,y,z$ whose coefficients are quadratic forms
in $\lambda$ and $\mu$. Moreover $D(\lambda,0)=\lambda^2C_1$, $E(0,\mu)=\mu^2C_2$, and
$D(\lambda_i,\mu_i)=E(\lambda_i,\mu_i)=C_i$ for $i=3,4$.
\end{cor}

In turns out that we can learn more about $D$ and $E$. The proof of the following proposition 
is taken from \cite{St}.

\begin{prop}
\label{stipins}
In Corollary \ref{form} we also have $D(0,1)=uC_2$ and $E(1,0)=vC_1$ for some $u,v\in\mathbb C$.
\end{prop}
\begin{proof}
We prove that the polynomial $D(0,1)$ is proportional to $E(0,1)$. Then the first formula will follow
from Corollary \ref{form} and the second one can be proved similarly.

Consider the pencil $G(a,b)=aD(0,1)+bE(0,1)$. It is projective 
(see \cite{Co}, p.32) to
the degenerate pencil generated by 0 and $F_2$. By the Chasles theorem (Ibid.) the
correspondence $G(a,b)\leftrightarrow aF_2$ is one-to-one which is  impossible
unless $D(0,1)$ and $E(0,1)$ are proportional.
\end{proof}
\begin{cor}
\label{eqn}

The polynomials $D$ and $E$ can be expressed in terms of $C_1$ and $C_2$ as follows.
\begin{eqnarray}
\label{1}
D(\lambda,\mu)=\lambda^2C_1+\lambda\mu X+u\mu^2 C_2 \\
E(\lambda,\mu)=v\lambda^2C_1+\lambda\mu Y+\mu^2C_2,
\end{eqnarray}
where $X$ and $Y$ are some forms in $x,y,z$ of degree $2d-6$ .
\end{cor}
Now we are ready to prove the strongest property of the pencil of the Hessians.

\begin{prop}
\label{const}
There exists a form $C$ of degree $2d-6$ in $x,y,z$ and quadratic forms $p$ and $q$ in $\lambda$ and $\mu$ such that
\begin{eqnarray}
\label{2}
H(\lambda,\mu)=C(p\lambda F_1+q\mu F_2)
\end{eqnarray}
for all $(\lambda,\mu)\in\mathbb C^2$.
\end{prop}
\begin{proof}
Assume without any loss that $C_1$ and $C_2$ are not identically zero.
Recalling that $D(\lambda_i,\mu_i)=E(\lambda_i.\mu_i)$ for $i=3,4$, and
using Corollary \ref{eqn} we obtain the following two equalities: 
\begin{eqnarray}
\label{3}
\lambda_i\mu_i(X-Y)=C_1\lambda_i^2(v-1)+C_2\mu_i^2(1-u)
\end{eqnarray}
where $i=3,4$.
Since $\lambda_i\mu_i\not=0$ for $i=3,4$ and also 
$\lambda_3\mu_4-\lambda_4\mu_3\not=0$ we can eliminate $Y-X$ and obtain 
\begin{eqnarray}
\label{4}
(v-1)\lambda_3\lambda_4C_1=(1-u)\mu_3\mu_4C_2.
\end{eqnarray}

The equalities $v=u=1$ would imply that $X=Y$ 
 whence $D=E$ and  $F(\lambda,\mu)$ divides
$H(\lambda,\mu)$ for all $[\lambda:\mu]$. Since this is impossible for generic fibers
\eqref{4} implies that $C_1$ and $C_2$ are proportional.
Since we used about $F_1$ and $F_2$ only that they are special we can change
the basis of the pencil and extend this statement to an arbitrary pair 
 $(C_i,C_j)$, $i\not=j$ . Choosing $C=C_1$ we see that every $C_i$ is proportional
 to $C$. 

Now consider again the equalities (1) and (2). Plugging in $\lambda=\lambda_3$
and $\mu=\mu_3$ yields (together with Corollary \ref{form}) that $X$ and $Y$ are proportional to $C$. Now factoring out $C$ implies \eqref{2}.
\end{proof}

\section{The main result}

In this section we proof the following main result of the note.

\begin{thm}
\label{mainer}
Let $\mathcal P$ be a pencil of plane curves of degree $d>1$ with irreducible generic fiber. Then

(i) if it has three completely reducible fibers then it can have at most one
more special fiber;

(ii) if it has three completely reducible fibers and another special fiber then
this special fiber is completely reducible and every two completely reducible
fibers intersect transversally (i.e., $|B|=d^2$). 
\end{thm}

\begin{proof}
(i) The first statement follows immediately
from Proposition \ref{const}. Indeed for each special fiber
its Hessian is divisible by the fiber. Thus the respective values of $\lambda$ and $\mu$ form a root of the polynomial
$$r(\lambda,\mu)=p(\lambda,\mu)\lambda\mu-q(\lambda,\mu)\lambda\mu$$
of degree four. If there are five such fibers then $r$ is identically zero
whence
$$H(\lambda,\mu)=p(\lambda,\mu)CF(\lambda,\mu)$$
which is impossible for generic $(\lambda,\mu)$.

(ii) Notice that the first part of the statement (ii) follows from the second 
one since non-reduced curves cannot have intersection with $d^2$ points.
Thus it suffices to prove the second part.

Now we assume that $\mathcal P$ has three completely reducible fibers
and another special fiber.  First suppose that $d=3$, i.e., $C=\emptyset$
 where $C$ is the curve from \eqref {2}.
 Then the equality \eqref{2} implies that for each cubic which is
 a fiber of the pencil its Hessian is also a fiber. 
On the other hand, it is well-known that the Hesse pencil has any plane
 cubic as a fiber. Thus taking a generic fiber of our pencil and its
Hessian we obtain two generators of our pencil and of the Hesse pencil whence
these two pencils coincide. The statement (ii) for the Hasse pencil is well-known.

Now suppose that $d>3$, i.e., $C\not=\emptyset$, and there exists
$P\in B$ at which the multiplicity $m$ of any fiber is greater than 1.
We assume as before $P=[0:0:1]$ and
let  $\alpha$ be a linear factor of $F^{(m)}$ (i.e., a tangent line to $F$ at $P$)
where $F$ is a generic
fiber of $\mathcal P$.
According to Lemma \ref{hessian}, $F^{(m)}$ divides $H(F)^{(3m-4)}$ whence so does 
 $\alpha$, i.e., $\alpha$ is a tangent line to $H(F)$ at $P$. 
  Since $F$ is generic we have $H(F)=CG$ where $G$ is a fiber
different from $F$. Since $\alpha$ is tangent to $H(F)$ but not to $G$
(by Lemma \ref{mult}) it is tangent to $C$. In other words we have proved
that every tangent line to every generic fiber at $P$ is tangent to $C$
whence there exist fibers with common tngent lines at $P$. This 
contradicts Lemma \ref{mult}.
\end{proof}

The following problem still remains unsolved.

{\bf Problem.}
Does there exist a pencil of a degree greater than 3 with four completely reducible fibers?

By the previous theorem the fibers of this pencil must intersect
transversally, i.e., its existence is equivalent
to the existence of a ($4,d)$-net in $\mathbb P^2$
with $d>3$.

The negative solution of this problem in \cite{St} has a gap.


\begin{thebibliography}{9}

\bibitem{Bo}
A. Bodin, \emph{Reducibility of rational functions in several variables}, math.AG/0510434.

\bibitem{Co}
Coolidge, J. L. A Treatise on Algebraic Plane Curves, New York, Dover, 1959.

\bibitem{Di}
A. Dimca, \emph{Pencils of Plane Curves and Characteristic
Varieties}, math.AG/0606442.


\bibitem{FY}
M. Falk and S. Yuzvinsky, \emph{Multinets, Ressonance Varieties, and
pencils of plane curves}, Compositio Math.  143 (2007), 1069-1088.


\bibitem{Fu}
W.Fulton, \emph{Algebraic curves}, Benjamin/Cummings, 1969.

\bibitem{Ha}
J. Hadamard, \emph{Sur les conditions de d\'ecomposition des formes},
Bull. SMF \textbf{27} (1899), 34-47.

\bibitem{Hal}
G. Halphen,\emph{Oeuvres de G.-H, Halphen}, t. III, Gauthier-Villars, 1921, 1-260.

\bibitem{LY}
A. Libgober and S. Yuzvinsky, \emph{Cohomology of the Orlik-Solomon algebras and local systems},
Compositio Math., 121 (2000), 337-361.



\bibitem{PY}
J. V. Pereira and S. Yuzvinsky, \emph{Completely reducible hypersurfaces in a pencil},
math.AG/0701312.

\bibitem{St}
J. Stipins, \emph{On finite $k$-nets in the complex projective plane}, Ph. D. thesis,
The University of Michigan, 2007.

\bibitem{Vi}
A. Vistoli, \emph{The number of reducible hypersurfaces in a pencil}, Invent. Math.
\textbf{112} (1993), 247-262.


\bibitem{Yu}
S. Yuzvinsky, \emph{Realization of finite Abelian groups by nets in
$\mathbb{P}^2$}, Compositio Math. \textbf{140} (2004),
1614--1624.


\end{thebibliography}
\end{document}